\documentclass[a4paper]{article}
\usepackage{RR}
\RRNo{6284}
\usepackage{hyperref}
\usepackage[]{amssymb}
\usepackage[]{amsmath}
\usepackage{psfrag}

\newtheorem{theorem}{Theorem}

\newtheorem{hypothesis}{Hypothesis}
\newtheorem{proposition}{Proposition}
 
\graphicspath{{images/}}

\newenvironment{proof}{\paragraph{Proof:}}{\endpf\\}
\def\endpf{\hfill$\Box$\medskip}

\RRdate{Septembre 2007}
\RRauthor{
Sapna Nundloll 
  \thanks[sfn]{INRIA - COMORE, 2004 route des Lucioles, BP93, F-06902 Sophia Antipolis Cedex, France.}%
\and 
Ludovic Mailleret 
\thanks{INRA - UR880, 400 route des Chappes, BP167, F-06903, Sophia Antipolis Cedex, France.}
\and
Fr\'ed\'eric Grognard
\thanksref{sfn}
}
\authorhead{S. Nundloll, L. Mailleret \& F. Grognard}
\RRtitle{Effets de la r\'ecolte partielle des cultures sur la lutte biologique}
\RRetitle{The Effects of Partial Crop Harvest on Biological Pest Control}
\titlehead{The effect of partial crop harvest on biological pest control}
\RRresume{Ce document montre les effets de la r\'ecolte partielle p\'eriodique sur la lutte biologique inondative dans le cas de la protection d'une plantation \`a croissance continue. Comme la r\'ecolte partielle est susceptible de pr\'elever une partie des ravageurs ainsi que des auxiliaires de lutte biologique, on ne peut ignorer son influence sur son efficacit\'e.

 Le mod\`ele \'etudi\'e consiste en un mod\`ele proie-pr\'edateur en EDO
 classique d\'ecrivant l'interaction biologique entre les deux populations,
 auquel s'ajoute une partie discr\`ete repr\'sentant les ph\'enom\`enes,
 intrins\`quement discrets, li\'es \`a la r\'ecolte et \`a la lutte biologique. 
 Les param\`etres de r\'ecolte partielle sont suppos\'es fix\'es par la
 physiologie de la plante ainsi que la demande \'economique; la p\'eriode de
 l\^acher des auxiliaires et le budget investi sont donc les seuls
 param\^etres modifiables du syst\`eme. L'exigence de stabilit\'e de l'\'etat du
 syst\`eme o\`u les ravageurs sont absents nous donne un budget minimal
 d'auxiliaires \`a utiliser par unit\'e de temps qui peut \^etre fonction de la
 p\'eriode des l\^achers. Nous d\'emontrons que la p\'eriode des l\^achers influe sur le budget minimal quand elle a lieu plus souvent que la r\'ecolte, mais n'a aucun effet quand elle a lieu moins fr\'equemment.
}

\RRabstract{In this paper, the effects of periodic partial harvesting of a continuously grown crop on augmentative biological control are analyzed. Partial harvesting can remove a proportion of both pests and biological control agents, so its influence on the control efficiency cannot be {\it a priori} neglected. An impulsive model consisting of a general predator-prey model in \textsc{ode}, augmented by a discrete component to depict releases of biological control agents and the periodic partial harvesting is used. The periods are taken as integer multiples of each other. A stability condition for pest eradication is expressed as the minimal value of the budget per unit time to spend on predators. We consider the partial harvesting period to be fixed by both the plant's physiology and market forces so that the only manipulated variable is the release period. It is shown that varying the release period with respect to the harvest period influences the minimal budget value when the former is carried out more often than the latter and has no effect when releases take place as often as or less frequently than the partial harvests.}
\RRmotcle{Dynamiques proie-pr\'edateur; r\'ecolte partielle; lutte biologique inondative; commande impulsive; stabilit\'e; budget minimal}
\RRkeyword{Predator-prey dynamics; partial harvest; inundative biological control; impulsive control; stability; minimal budget}
\RRprojets{Comore}
\RRtheme{\THBio} 
\URSophia 
\begin{document}
\makeRR   
\section{Introduction}

Biological control is the reduction of pest populations to harmless levels through the release of their natural enemies. The latter can include both parasitic and predatory species, which are deployed at selected locations throughout the crop and, wherever possible, to specific parts of individual plants where the pest is likely to attack. Successful control projects in the field have involved the use of only one predatory species such as in \cite{EHL05,JACCROFEN01}, as well as more complex biodiverse schemes such as those suggested by \cite{MUReA85,JACeA01,EHL90,ROSeA95} and the references therein. The target pest species and the setting, i.e. where the crop is grown, usually determines the type of control required, namely whether pest eradication is necessary or not. For an exhaustive list of definitions and applications, we refer the reader to \cite{EILHAJLOM01,RODNAV03}.

In this report, we consider the protection of continuously grown crops which have zero tolerance to pest invasions. There are two aspects in this type of culture. 

\begin{itemize}
\item Firstly, inundative control which is a prophylactic method of pest control yields the most satisfatory results when implemented (see \cite{FENeA01,SKIeA02,FENeA02}  for theoretical/simulatory studies and \cite{EHLeA97,JACeA01,JACCROFEN01,DCW01} for real life experiments). A calculated number of predators are repeatedly injected into the ecosystem, independently of the detection of pest insects in the greenhouse. Such populations are not allowed to thrive and consist only of individuals whose main source of subsistence is the pest insect, in the absence of which, they (the predator insects) rapidly die out.  The frequency of the releases and the number of predators injected each time ensures that a minimal 'sentry' population is present to reduce the damage caused by the pests on their attack. 
\item Secondly, over their growing period, these crops are partially harvested on a regular basis. Since it is known that harvests are likely to influence, even counterintuitively, predator-prey dynamics  \cite{SCU78,NEGGAK07}, it has to be taken into account in the formulation of the problem.
\end{itemize}

We consider the simplest ditrophic case whereby one predatory species is used to eradicate a pest population. Our model consists of ODEs augmented by a discrete component to incorporate the effect of partial harvest and releases that by their very nature are discrete phenomena. This is a classical formulation that is used widely in the literature where impulsive dynamics are studied. Examples are \cite{LIUZHACHE05a,NEGGAK07} in the context of agricultural ecosystems, \cite{SHUSTOAGU98} in epidemiology, \cite{LAKARI00} in pulsed chemotherapy to cite some. Few papers in the literature on impulsive crop protection however seem to focus on stability of the pest-free state: yet this is of practical importance especially for high valued crop cultures. 

In our work, we attempt to give an economic dimension to the solution of our problem by defining the releases in terms of the number of predators to invest in over a budget period.
Using Floquet Theory as presented by \cite{SHUSTOAGU98}, we are able to express the stability condition as the minimal number of predators per budget period required to drive the pests to zero at a given release frequency. \cite{MAIGRO06} showed how this number varied with the release period chosen. The worst case scenario of pest attack occurring at an intermediate stage between two predator releases was considered and the optimal release policy which would guarantee the most efficient protection against surges in the pest population was calculated. In particular, it is shown that the higher the frequencies of predator release, the smaller the time interval over which the pest population was above a threshold commonly referred to as the Economic Injury Level \cite{STEeA59} - and hence the lower the damage incurred by the crops. 

In line with the work of \cite{MAIGRO06}, we investigate how the frequency of releases is to be varied with respect to the (fixed) harvesting frequency to minimise the minimal budget value. We consider the harvest period as a reference since it is set by market constraints. The effect of partial harvesting is similar to that of pesticide usage proposed by \cite{LIUZHACHE05a} in their Integrated Pest Management (IPM) strategy. Our model departs from the latter's in three ways. Firstly, both the predator and pest populations are subjected to partial harvesting when this occurs. Secondly, hypotheses made on the functions governing the population changes are weak and can encompass most of the density-dependent functions proposed in the literature. Finally, one period is taken as the integer multiple of other. This feature is key in solving for the stability condition to obtain the minimal budget value. The case where the frequencies are not the same is included. 

It is shown that for a given harvest period, when releases take place less often or as often as harvests, the minimal budget is at a calculated value which is independent of release period. However, when releases take place more often than harvests, the minimal budget required always exceeds this value. This result runs counter with that obtained by \cite{MAIGRO06}: merging the two seems to indicate that the harvest frequency is a threshold that should not be exceeded when releasing predators for efficient biological control.

In the first section of this article, the system model is presented. The mathematical analysis of the system's stability and the formulation of the stability condition in terms of the minimal budget are presented in the next section. A brief interpretation of the mathematical results follows. Finally, we conclude with a discussion on their implications.

\section{Model description}
The model we present consists of a continuous part to depict the predator-prey interaction. We consider the case at the onset of pest invasion where the crop - the pest food supply - is in abundance. Because of this, at this stage, it is sufficient to model only the pest $x$ and predator $y$ species. 

\begin{equation}
\label{model}
\left\{
	\begin{array}{cclc}
		\dot{x} & = & f(x) - g(x)y & \\
		\dot{y} & = & h(x)y - dy & \\
		x(nT_h^+) & = & (1-\alpha_x)x(nT_h) & \forall n \in\mathbb{N}\\
		y(nT_h^+) & = & (1-\alpha_y)y(nT_h)+\delta\left(nT_h\mbox{ mod }T_r \right)\mu T_r &\forall n\,\in\mathbb{N}\\
		y(mT_r^+) & = & (1-\delta\left(mT_r\mbox{ mod }T_h \right)\alpha_y)y(mT_r)+\mu T_r &\forall m\,\in\mathbb{N}
	\end{array}\\
\right.
\end{equation}

The first two equations govern the intrinsic predator-prey interaction occurring in the system. The three ones depict the impulsive phenomena that we consider with harvest taking place at $nT_h$ and releases at $mT_r$. 

In the continuous part, the functions discussed are not specified so they are representative of as many systems as possible. Only the following hypotheses are made. 

\begin{hypothesis}
	\label{hyp1}
$~$
Let $f(x)$, $g(x)$ and $h(x)$ be locally Lipschitz continuous in
$\mathbb{R^+}$ such that 
	\begin{itemize}
		\item $f(0) = 0$
		\item $g(0) = 0$, $g'(0) > 0 $ and $g(x) > 0$ $\forall x>0 $
		\item $h(0) = 0$ and $h(x) > 0$ $\forall x>0$
		\item $\frac{f(x)}{g(x)}$ and $\frac{g(x)}{x}$ are upper bounded for $x\geq 0$
	\end{itemize}
\end{hypothesis}

$f(x)$ is 
 the growth velocity or feeding input of the pests. It represents the growth function of the pest species and in our model, it also encompasses any non-predatory losses of the pest population (e.g. logistic growth). We assume that the predator population is never large enough for intra-predator interaction to take place so the functional and numerical responses can be expressed solely in terms of the prey numbers, i.e. as $g(x)$ and $h(x)$ respectively.

We assume that pest growth rate, the functional and numerical responses are all nil when the ecosystem is pest-free. 

The functional response is increasing for small pest population levels. We also
consider that, in the presence of pests, predation always takes place with a
negative impact on $x$ ($g(x)>0$) and a positive impact on $y$
 ($h(x)>0$). Note that conditions can be induced as much by the predator
 insect foraging abilities \textit{per se} as they can be facilitated by
 placing the predator insects at known locations on the plant where the pests
 are most likely to attack. In classical density dependent models, $g(x)$ is
 bounded or linear, so that $\frac{g(x)}{x}$ is always bounded. The
 boundedness of $\frac{f(x)}{g(x)}$ means that there is no value of $x$ where
 the pest growth $f(x)$ overwhelmingly dominates the predation $g(x)$, which
 would render the biological control impossible. 

Partial crop harvests and predator releases occur respectively every $T_h$ and $T_r$. $\alpha_x$ and $\alpha_y$ represent the respective proportions of the prey and predator populations affected at each harvest. These parameters are allowed be different since in reality, it is very likely that each species tends to occupy different parts of the plant. We also assume that the insects are uniformly distributed throughout our plantation so that the effect of partial harvesting is directly correlated with the number of plants harvested. We assume linear maturation of the crop so the proportion of crops harvested each time and hence insects removed is considered as fixed. The $\delta$-function is defined thus to identify instants of simultaneous partial harvest and predator release.

\begin{equation}
	\begin{array}{lll}
		\delta(\theta) & = &
		\left\{
			\begin{array}{ll}
 				1 & \mbox{ if } \quad \theta = 0\\
 				0 & \mbox{ otherwise }
			\end{array}
		\right.
	\end{array}
\end{equation} 

Finally, we presume that we have a fixed budget of predators over a designated time period that is distributed evenly among the releases that are carried out. $\mu$ refers to the total number of predators purchased per time unit. Expressing $T_r$ in the same units as the budget period gives the control $\mu T_r$ as the number of predators released every $T_r$. 

\section{Mathematical analysis}

In our analysis, we restrict ourselves to the case where either one of the periods (release or partial harvests) is the integer multiple of the other. Note however that the model (\ref{model}) formalism is more general. We study the system in the absence of pests, i.e. when $x = 0$. In addition of being invariant, it is the target state of our system. The stability of the system around that state is therefore of interest. Our analysis takes place separately for the case when releases are more frequent than harvests, and when they are less frequent.

We show that in the absence of pests at the initial time, the predator population converges towards a positive periodic solution. We then demonstrate that when preys are present at the initial time, convergence of the predator population also takes place to that same periodic solution, while the preys go extinct provided some condition on the parameters is verified.

\subsection{Pest-free stability analysis}
\subsubsection*{Releases more frequent than harvests}

\begin{proposition}
\label{p_Tr}
Let $T_h = kT_r$ where $k \in \mathbb{N}^*$ and Hypotheses \ref{hyp1} be satisfied. Then, in the absence of pests, model (\ref{model}) possesses a globally stable periodic solution 

\begin{equation}
\label{theo_periodicsol_hi}
\left(
	x_{ph}	\left(
			t
		\right)
	, y_{ph}	\left(
			t
		\right)
\right)
= 
\left(
	0, y^*e^{-d (t\mod T_h)}+\mu T_r e^{-d (t\mod T_r)}\sum_{j = 0}^{\lfloor{ \frac{t\mod T_h}{T_r}}\rfloor  -
  1}{e^{-jdT_r}}
\right)
\end{equation}

where

\begin{equation}
\label{lim_y_hi}
y^* = \frac{	
			\left(
				\frac{1 - e^{-dT_h}}{1 - e^{-dT_r}} 
			\right) (1 - \alpha_y) + \alpha_y }{1 -(1 - \alpha_y)e^{-dT_h} }
		\mu T_r
\end{equation}

\end{proposition}

\begin{proof}

When $T_h = kT_r$, in the absence of pests and using Hypotheses \ref{hyp1}, the system is simplified to

\begin{equation}
\label{model_hi}
\left\{
	\begin{array}{ccl}
		\dot{x} & = & 0 \\
		\dot{y} & = & -dy \\
		x(mT_r^+) & = & (1- \delta\left(m \mbox{ mod }k \right) \alpha_x)x(mT_r) \\
		y(mT_r^+) & = & (1- \delta\left(m \mbox{ mod }k \right) \alpha_y)y(mT_r) + \mu T_r \\
		&  &\forall m\,\in\mathbb{N}
	\end{array}\\
\right.
\end{equation}

The pest population stays nil since in the absence of pests, their population does not change either. The solution 
$$
x_{ph}(t) = 0
$$
is trivial.

On the other hand, the predator population will vary according to the number of predators manually injected into the system and, since the population is non-zero, according to the partial harvest effect. The absence of their source of food will cause an exponential decay of the population. We demonstrate that these forces will provoke the predator population to reach a periodic pattern of period equal to $T_h$, which we shall refer to as the \textit{period of reference}. The instant following a coinciding partial harvest and release is taken as the \textit{point of reference}.

To prove Proposition \ref{p_Tr}, we first show by induction that the predator population right after a release can be expressed in terms of the point of reference as follows
\begin{equation}
\label{y_Tr}
\begin{array}{ccl}
y(nT_h + iT_r^+) & = & y(nT_h^+)e^{-idT_r} + \mu T_r \displaystyle\sum_{j = 0}^{i - 1}{e^{-jdT_r}}
\end{array}
\end{equation}
where $i \in [0, 1, \hdots, (k - 1)]$

It is seen that (\ref{y_Tr}) is valid for $i=0$ since it is equal to
\begin{equation*}
\label{y_rec_init_hi}
\begin{array}{ccl}
  y(nT_h^+) & = & y(nT_h^+)e^{0} + \mu T_r \displaystyle\sum_{j = 0}^{ - 1}{e^{-jdT_r}} = y(nT_h^+)\\
\end{array}
\end{equation*}
Now suppose that (\ref{y_Tr}) holds for $i = q$ where $q \in [0, 1, \hdots, k-2 ]$, i.e.
\begin{equation*}
\label{y_rec_q_hi}
\begin{array}{ccl}
y(nT_h + qT_r^+) & = & y(nT_h^+)e^{-qdT_r} + + \mu T_r \displaystyle\sum_{j = 0}^{q - 1}{e^{-jdT_r}}
\end{array}
\end{equation*}
We will now show that (\ref{y_Tr}) is valid for $i = q + 1$. We calculate $y(nT_h + (q + 1)T_r^+) $ from $y(nT_h + qT_r^+)$ using (\ref{model_hi}), then substituting (\ref{y_rec_q_hi}) into (\ref{y_rec_qplus1_hi}) as follows
\begin{equation}
 \label{y_rec_qplus1_hi}
\begin{array}{ccl}
 y(nT_h + (q + 1)T_r^+) & = & y(nT_h + qT_r^+)e^{-dT_r} + \mu T_r\\
			& = & \left( y(nT_h^+)e^{-qdT_r} + \mu T_r \displaystyle\sum_{j = 0}^{q - 1}{e^{-jdT_r}}\right)e^{-dT_r} + \mu T_r\\
			& = & y(nT_h^+)e^{-(q + 1)dT_r} + \mu T_r \displaystyle\sum_{j = 1}^{q}{e^{-jdT_r}} + \mu T_r\\
			& = & y(nT_h^+)e^{-(q + 1)dT_r} + \mu T_r \displaystyle\sum_{j = 0}^{q}{e^{-jdT_r}} 
\end{array}
\end{equation}
so that (\ref{y_Tr}) holds true for $i \in [0, 1, \hdots, k-1] $.

To evaluate the evolution of $y$ according to the period of reference $T_h$, we need to calculate the value of $y((n + 1)T_h^+) $, which is equivalent to $y(nT_h + kT_r^+)$, in terms of $y(nT_h^+)$ . At this point however, partial harvesting takes place before predator release; so we first express  it in terms of $y(nT_h + (k - 1)T_r^+)$ then expand the expression using (\ref{y_Tr}) as follows

\begin{equation*}
\begin{array}{ccl}
 y\left( (n + 1) T_h^+ \right) 	& = & y\left( nT_h + (k - 1)T_r^+ \right)e^{-dT_r}(1 - \alpha_y) + \mu T_r\\
				& = & \left( y(nT_h^+)e^{-d(k - 1)T_r} + \mu T_r \displaystyle\sum_{j = 0}^{k - 2}{e^{-jdT_r}} \right)e^{-dT_r}(1 - \alpha_y) + \mu T_r\\
				& = & y(nT_h^+)e^{-dT_h} + \mu T_r (1 - \alpha_y)\displaystyle\sum_{j = 1}^{k - 1}{e^{-jdT_r} + \mu T_r}\\
				& = & y(nT_h^+)e^{-dT_h} + \mu T_r \left( (1 - \alpha_y)\displaystyle\sum_{j = 0}^{k - 1}{e^{-jdT_r}} + \alpha_y \right)
\end{array}
\end{equation*}
Note that the summation term can also be evaluated so the sequence is expressible as
\begin{equation}
\label{y_series_hi}
\begin{array}{ccl}
 y\left( (n + 1)T_h^+ \right)	& = & 	y(nT_h^+)e^{-dT_h} + \mu T_r 
								\left((1 - \alpha_y	)\displaystyle\frac{1 - 		e^{-dT_h}}{1 - e^{-dT_r}} + 		\alpha_y 	  	
								\right)
\end{array}
\end{equation}
In this linear dynamical system, the coefficient of $y(nT_h^+)$, $e^{-dT_h} $ is less than one in magnitude, so the sequence will converge to a limit, the equilibrium of (\ref{y_series_hi}). This equilibrium yields (\ref{lim_y_hi}) and the convergence of $y(t)$ to a periodic solution $y_{ph}(t)$ based on $y^*$.

Now that we have established the existence of the periodic solution $y_{ph}(t)$,
we seek to formulate it. We focus on a reference period over $nT_h < t \leq (n
+ 1)T_h$ during which  $y_{ph}(t)$ is piecewise continuous, with the continuous components separated by predator releases. The continuous intervals are defined over $nT_h + i T_r < t \leq nT_h + (i+1)T_r$ where $i \in [0, 1,\hdots, k-1]$. For a given value of $t$, the value of $i$ is easily identified
as being $i=\lfloor{ \frac{t\mod T_h}{T_r}}\rfloor$. The value of $y_{ph}(t)$
is then of the form
\[
y_{ph}(t)=y_{ph}(nT_h+iT_r^+)e^{-d (t\mod T_r)}
\]
and, from (\ref{y_Tr}) with $y(nT_h^+)=y^*$, we have that
\[
y_{ph}(nT_h+iT_r^+)=y^*e^{-idT_r} + \mu T_r \displaystyle\sum_{j = 0}^{i - 1}{e^{-jdT_r}}
\]
so that
\[
\begin{array}{lll}
y_{ph}(t)&=&\left(y^*e^{-idT_r} + \mu T_r \displaystyle\sum_{j = 0}^{i -
  1}{e^{-jdT_r}}\right)e^{-d (t\mod T_r)}\\
&=&y^*e^{-d (t\mod T_h)}+\mu T_r e^{-d (t\mod T_r)}\displaystyle\sum_{j = 0}^{i -
  1}{e^{-jdT_r}}\\
&=&y^*e^{-d (t\mod T_h)}+\mu T_r e^{-d (t\mod T_r)}\displaystyle\sum_{j = 0}^{\lfloor{ \frac{t\mod T_h}{T_r}}\rfloor  -
  1}{e^{-jdT_r}}
\end{array}
\]
This is of the same form as proposed in (\ref{theo_periodicsol_hi}), thereby completing our proof. 
\end{proof}

\begin{figure}[htbp]
\begin{center}
\includegraphics[height=9cm]{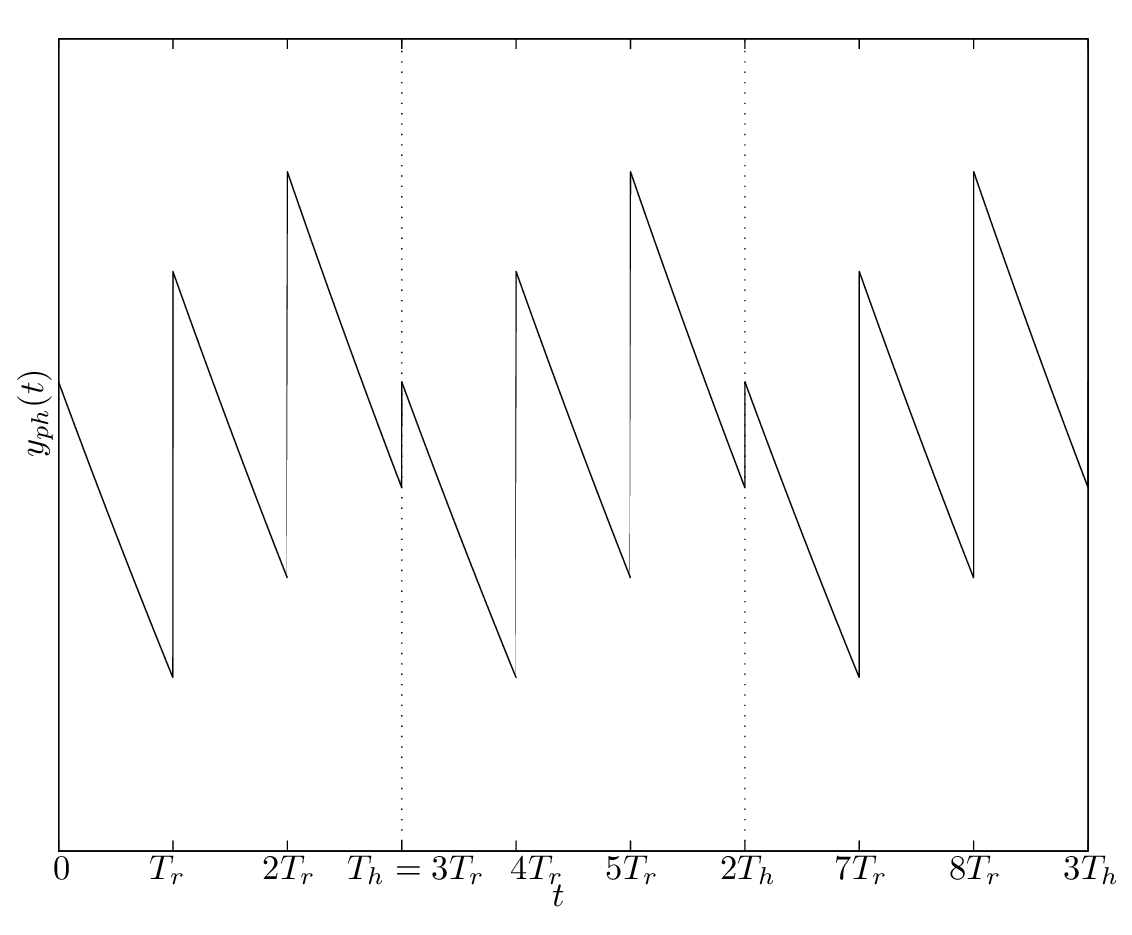}
\caption{Form of the periodic solution $y_{ph}(t)$ in the case where $k=3$. Releases of predators are apparent at every $mT_r$ instant,
  while the cumulative effect of harvest and release leads to an apparent
  smaller release at every $nT_h$ instant.
Between those instants, the population decays exponentially since it has no
prey to feed on.}\label{fig:yph}
\end{center}
\end{figure}

The form of the $y_{ph}$ function is illustrated on Figure \ref{fig:yph}.  


\subsubsection*{Releases less frequent than harvests}

When harvesting is more frequent than the release of predators, we have a similar result about the existence of a periodic solution.

\begin{proposition}
Let $T_r = kT_h$ where $k \in \mathbb{N^*}$ and Hypotheses \ref{hyp1} be satisfied.
Then, in the absence of pests, model (\ref{model}) possesses a globally stable periodic solution 
\begin{equation}
\label{theo_periodicsol_lo}
\left(
	x_{pr}	\left(
			t
		\right)
	, y_{pr}	\left(
			t
		\right)
\right)
= 
\left(
	0, y^*e^{-d(t\mod T_r)}(1 - \alpha_y)^{\lfloor
							\frac{t\mod T_r}{T_h}
						\rfloor
						}
\right)
\end{equation}
where 
\begin{equation}
\label{lim_y_lo}
	y^* = \frac{\mu T_r}{1 - (1 - \alpha_y)^ke^{-dT_r}}
\end{equation}

\end{proposition}

\begin{proof}

When $T_r = kT_h$, in the absence of pests and using Hypotheses \ref{hyp1}, the system is simplified to
\begin{equation}
\label{model_lo}
\left\{
	\begin{array}{ccl}
		\dot{x} & = & 0 \\
		\dot{y} & = & -dy \\
		x(nT_h^+) & = & (1- \alpha_x)x(nT_h) \\
		y(nT_h^+) & = & (1- \alpha_y)y(nT_h)+ \mu T_r\delta\left(n \mbox{ mod }k \right) \\
		&  &\forall n\,\in\mathbb{N}
	\end{array}\\
\right.
\end{equation}

As previously explained, $x_{pr}(t)$ is solved for trivially as being 

$$x_{pr}(t) = 0$$


We prove that the predator population again reaches a periodic solution. This time, however, the \textit{period of reference} is $T_r$. The \textit{point of reference} is the instant after a coinciding harvest and release. We show by induction that the population after a harvest can be expressed as
\begin{equation}
\label{y_Th}
y(mT_r + iT_h^+) = y(mT_r^+)e^{-idT_h}(1 - \alpha_y)^i
\end{equation}  
where $i \in [0, 1, \hdots, (k - 1)]$.  

It is seen that (\ref{y_Th}) is valid for $i= 0$ since it resumes to
\begin{equation*}
\label{y_rec_init_lo}
\begin{array}{ccl}
y(mT_r^+) 	& = & y(mT_r^+)e^0(1 - \alpha_y)^0\\
\end{array}
\end{equation*}
Suppose (\ref{y_Th}) holds for $i = q$ where $q \in [0, 1, \hdots k-2]$, i.e. 
\begin{equation}
\label{y_rec_q_lo}
\begin{array}{ccl}
y(mT_r + qT_h^+) & = & y(mT_r^+)e^{-qdT_h}(1 - \alpha_y)^q
\end{array}
\end{equation}
We will now show that (\ref{y_Th}) is valid for $i = q + 1$. We calculate the value of $y$ when $i = q + 1$ in terms of $y(nT_h + qT_r^+)$, knowing from $\dot{y} = -dy$ in (\ref{model_lo}) it will be an exponential decay with the added component for the harvest. We then substitute (\ref{y_rec_q_lo}) in and obtain
\begin{equation}
\label{y_rec_qplus1_lo}
\begin{array}{ccl}
y(mT_r + (q + 1)T_h^+) & = & y(mT_r + qT_h^+)e^{-dT_h}(1 - \alpha_y)\\
			& = & \left( y(mT_r^+)e^{-qdT_h}(1 - \alpha_y)^q \right)e^{-dT_h}(1 - \alpha_y)\\
			& = & y(nT_h^+)e^{-d(q + 1)T_h}(1 - \alpha_y)^{q + 1}
\end{array}
\end{equation}
This is clearly the same form given from the expression in (\ref{y_Th}), thereby validating it.

$y^*$ is given as the fixed point of the sequence representing post-release instants. Therefore, using (\ref{y_Th}) for $i = k$ and model (\ref{model_lo}), we next calculate $y((m + 1)T_r^+)$ as
\begin{equation}
\label{y_series_lo}
\begin{array}{ccl}
y((m + 1)T_r^+) & = & y(mT_r + kT_h^+)\\
		& = & y(mT_r^+)e^{-kdT_h}(1 - \alpha_y)^k + \mu T_r\\
		& = & y(mT_r^+)e^{-dT_r}(1 - \alpha_y)^k + \mu T_r
\end{array}
\end{equation}

In this linear dynamical system, the coefficient of $y(mT_r^+)$,  $e^{-dT_r}(1 - \alpha_y)^k$ is less than one in magnitude, which confirms the existence of the fixed point $y^*$ to which the sequence converges. This equilibrium yields (\ref{lim_y_lo}) and the convergence of $y(t)$ to a periodic solution $y_{pr}(t)$.

Now that we have established the existence of the periodic solution $y_{pr}(t)$,
we seek to formulate it. We focus on a reference period over $mT_r < t \leq (m
+ 1)T_r$ during which  $y_{pr}(t)$ is piecewise continuous, with the
continuous components separated by harvests. The intervals of continuity
span $mT_r + i T_h < t \leq mT_r + (i+1)T_h$ where $i \in [0, 1,
\hdots, k-1]$. For a given value of $t$, the value of $i$ is easily identified
as being $i=\lfloor{ \frac{t\mod T_r}{T_h}}\rfloor$. The value of $y_{pr}(t)$
is then of the form
\[
y_{pr}(t)=y_{pr}(mT_r+iT_h^+)e^{-d (t\mod T_h)}
\]
and, from (\ref{y_Th}) with $y(mT_r^+)=y^*$, we have that
\[
y_{pr}(mT_r + iT_h^+) = y^*e^{-idT_h}(1 - \alpha_y)^i
\]
so that
\[
\begin{array}{lll}
y_{pr}(t)&=&\left( y^*e^{-idT_h}(1 - \alpha_y)^i \right)e^{-d (t\mod T_h)}\\
&=&y^*e^{-d (t\mod T_r)}(1 - \alpha_y)^{\lfloor{ \frac{t\mod
      T_r}{T_h}}\rfloor}
\end{array}
\]
which is exactly the expression given in  (\ref{theo_periodicsol_lo}) and
completes the proof.
\end{proof}

\begin{figure}[!ht]
\begin{center}
\includegraphics[height=9cm]{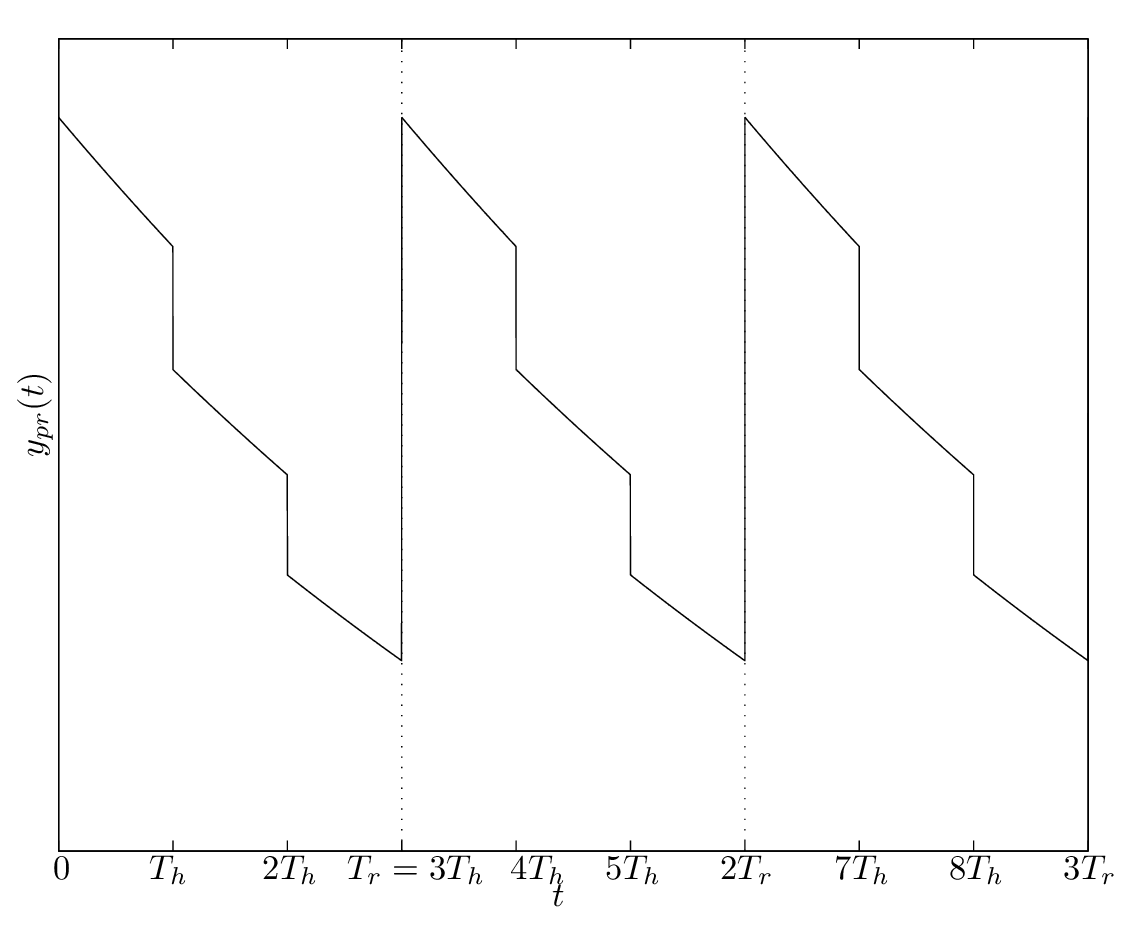}
\caption{Form of the periodic solution $y_{pr}(t)$ in the case where $k=3$. Harvests are apparent at every $nT_h$ instant,
  while the release of predators dominates the harvest at every $mT_r$ instant.
Between those instants, the population decays exponentially since it has no
prey to feed on.}\label{fig:ypr}
\end{center}
\end{figure}

The form of the $y_{pr}$ function is illustrated on Figure \ref{fig:ypr}.

\subsection{Global stability analysis}

Since we will study the convergence of the solutions to $(0,y_p(t))$ (where the $p$ subscript stands as well for $ph$ or $pr$), it will be convenient to describe the system in terms of the deviation coordinates
with respect to the reference periodic solution:
\begin{eqnarray*}
	\tilde{x}(t) &=& x(t) - x_p(t)\\
	\tilde{y}(t) &=& y(t) - y_p(t)
\end{eqnarray*}
This yields
\begin{eqnarray}
\dot{\tilde{x}} &=& f(x) - g(x)y\nonumber\\
		&=& f(\tilde{x}) - g(\tilde{x})(\tilde{y} + y_p(t))\label{xdot}
\end{eqnarray}
and
\begin{eqnarray}
\dot{\tilde{y}}	&=& h(x)y - dy - h(x_p)y_p + dy_p\nonumber\\
		&=& h(\tilde{x})(\tilde{y} + y_p(t)) - d\tilde{y}\label{ydot}
\end{eqnarray}

The impulsive effects on $\tilde x$ are obviously unchanged compared to those
on $x$. On the other hand, the release effects on $y$ disappear in $\tilde y$;
indeed, we have
\[
\tilde y(mT_r^+)=y(mT_r^+)-y_p(mT_r^+)=y(mT_r)+\mu T_r-(y_p(mT_r)+\mu T_r)=\tilde y(mT_r)
\]

The harvesting impulses are preserved in the expression of $\tilde y$
\[
\tilde y(nT_h^+)=y(nT_h^+)-y_p(nT_h^+)=(1-\alpha_y)y(nT_h)-(1-\alpha_y)y_p(nT_h)=(1-\alpha_y)\tilde y(nT_h)
\]

In the sequel, we will perform a global and a local stability analysis. For
the latter, we will need the computation of the linear approximation of the
deviation system around the periodic solution $(0,y_p(t))$:

\begin{equation}\label{linear}
\left\{
\begin{array}{lll}
\dot{\tilde{x}}&=&(f'(0)-g'(0)y_p(t))\tilde{x}\\
\dot{\tilde{y}}&=&h'(0)y_p(t)\tilde{x}-d\tilde{y}
\end{array}
\right.
\end{equation}

\subsubsection*{Releases more frequent than harvests}

We will first prove our result in the case where releases take place more
often than harvests. We obtain two different constraints for the Local
Asymptotic Stability (LAS) and Global Asymptotic Stability (GAS) of the
periodic solution in system
(\ref{model}). The latter is obviously stronger than the former, but is
sufficient in the case where pests outbreaks do not immediately take large
proportions.

In order to state the following theorem, we first need to define the 
function
\[
\underline{\mu}_h(S,r)=
	d\left(S+ \frac{\ln{\left(1
                  -\alpha_x\right)}}{rT_h}\right)\frac{1}{1- \left(\frac{\alpha_y\left(1 -
              e^{-dT_h}\right)}{1 - \left(1 -
              \alpha_y\right)e^{-dT_h}  }\right)\left(\frac{ e^{-\frac{d T_h}{k}}}{k\left(1 - e^{-\frac{dT_h }{k}}\right)  }\right)}
\] 
This function is increasing in $S$ and $r$ because the sign of the partial
derivatives is determined by the sign of the last factor, which can be shown
to be positive. Indeed, this factor is positive when

\[
\left(\frac{\alpha_y\left(1 -
              e^{-dT_h}\right)}{1 - \left(1 -
              \alpha_y\right)e^{-dT_h}  }\right)\left(\frac{ e^{-\frac{d T_h}{k}}}{k\left(1 - e^{-\frac{dT_h }{k}}\right)}\right)<1
\]  
and we have $\alpha_y e^{-\frac{d T_h}{k}}\leq \alpha_y$ and $1 - \left(1 -
              \alpha_y\right)e^{-dT_h}>\alpha_y$, so that
\[
\frac{\alpha_y e^{-\frac{d T_h}{k}}}{1 - \left(1 -
              \alpha_y\right)e^{-dT_h}}<1
\]
Also, $k\left(1 - e^{-\frac{dT_h }{k}}\right)\geq \left(1 -e^{-dT_h}\right)$
    since both sides of the inequality have the same value in $T_h=0$, and 
\[
\frac{d}{dT_h}\left(k\left(1 - e^{-\frac{dT_h
      }{k}}\right)\right)=de^{-\frac{dT_h }{k}}\geq de^{-dT_h }=\frac{d}{dT_h}\left(1 - e^{-dT_h
      }\right)
\]
which shows that $\underline{\mu}_h(S,r)$ is increasing in $S$ and $r$.

\begin{theorem}\label{theoremTh}

When $T_h=kT_r$ with $k\,\in\mathbb{N^*}$, the solution $(x(t),y(t))=(0,y_{ph}(t))$
of (\ref{model}) is LAS iff
\begin{equation}
\label{condThL}
  \mu > \underline{\mu}_h\left(\frac{f'(0)}{g'(0)},g'(0)\right)  
\end{equation}
and is GAS if
\begin{equation}\label{condTh}
  \mu > \underline{\mu}_h\left(\sup_{x\geq
  0}\frac{f(x)}{g(x)},\sup_{x\geq
  0}\frac{g(x)}{x}   \right)  
\end{equation}
\end{theorem}

\begin{proof}

We start with the proof of global convergence under condition (\ref{condTh}). In this proof, we will first show that $\tilde x$ goes to zero, from which we
will derive that $\tilde y$ goes to $0$ also (so that $y(t)$ converges to $y_{ph}(t)$).

Let the initial condition for system (\ref{xdot})-(\ref{ydot}) be $(\tilde
x_0,\tilde y_0)$ at time $t_0=0^+$, that is after the harvest and the predator
release that take place at the initial time. Analyzing (\ref{ydot}) and noting that $y_{ph}(t)+\tilde{y}=y(t)\geq 0$, we have
\[
\dot{\tilde{y}}\geq -d \tilde{y}
\]
so that $\tilde y(t)\geq \min(0,\tilde y_0) e^{-dt}$. 

In order to analyze the $\dot{\tilde x}$ equation, we define the function 
\begin{equation}\label{Gdef}
G(\tilde x)=\int_{x_0}^{\tilde x} \frac{1}{g(s)}ds
\end{equation}
which can easily be seen to be an increasing function of $\tilde x$ since
$g(s)>0$. Since we also have that \linebreak $g(s)<\left(\sup_{x\geq
  0}\frac{g(x)}{x}\right)s$, it is straightforward that
$\displaystyle\lim_{\tilde x\stackrel{>}{\rightarrow}0}G(\tilde x)=-\infty$. In order to show the
extinction of the pests we will then prove that $G(\tilde x)$ goes to $-\infty$ as $t$ goes to infinity. 
Therefore, we write the $G$ dynamics:
\[
\begin{array}{lll}
\frac{dG(\tilde x)}{dt}&=&\frac{1}{g(\tilde{x})}\dot{\tilde{x}}\\
&=&\frac{f(\tilde x)}{g(\tilde x) }-\tilde{y} - y_{ph}(t)\\
&\leq&\frac{f(\tilde x)}{g(\tilde x) }- \min(0,\tilde y_0) e^{-dt} - y_{ph}(t)
\end{array}
\]
We will now consider the evolution of $G$ between two successive harvests,
that is the evolution of $G$ between the times $nT_h^+$ and $(n+1)T_h$ for a
given $n$:
\[
G(\tilde{x}((n+1)T_h))\leq G(\tilde{x}(nT_h^+))+\int_{nT_h^+}^{(n+1)T_h}\left[
  \frac{f(\tilde{x}(s))}{g(\tilde{x}(s))}- \min(0,\tilde y_0) e^{-ds} -y_{ph}(s) \right]ds
\]
Since no impulse is present inside the integral, we can drop the $^+$
superscript in its lower extremity.

We will now analyze how the harvest that takes place at time $(n+1)T_h$
impacts $G$. We have
\begin{equation}\label{bound1}
\begin{array}{lll}
G(\tilde{x}((n+1)T_h^+))&=&\int_{{x}_0}^{\tilde{x}((n+1)T_h^+)} \frac{1}{g(s)}ds\\[.1cm]
&=&\int_{x_0}^{\tilde{x}((n+1)T_h)} \frac{1}{g(s)}ds+\int_{\tilde{x}((n+1)T_h)
}^{\tilde{x}((n+1)T_h^+)} \frac{1}{g(s)}ds\\[.2cm]
&\leq &G(\tilde{x}(nT_h^+))+\int_{nT_h}^{(n+1)T_h}\left[
  \frac{f(\tilde{x}(s))}{g(\tilde{x}(s))}- \min(0,\tilde y_0) e^{-ds} -y_{ph}(s) \right]ds \\
&&+\int_{\tilde{x}((n+1)T_h)
}^{\tilde{x}((n+1)T_h^+)} \frac{1}{g(s)}ds
\end{array}
\end{equation}

The last term represents the influence of harvest on $G$ and can easily be
approximated because \linebreak
$\tilde{x}((n+1)T_h)>\tilde{x}((n+1)T_h^+)=(1-\alpha_x)\tilde{x}((n+1)T_h)$. Denoting $S_g=\sup_{x\geq
  0}\frac{f(x)}{g(x)}$ and $r_g=\sup_{x\geq
  0}\frac{g(x)}{x}$, we have
\begin{equation}\label{bound2}
\int_{\tilde{x}((n+1)T_h)
}^{(1-\alpha_x)\tilde{x}((n+1)T_h)} \frac{1}{g(s)}ds\leq \int_{\tilde{x}((n+1)T_h)
}^{(1-\alpha_x)\tilde{x}((n+1)T_h)} \frac{1}{r_gs}ds=\frac{\ln(1-\alpha_x)}{r_g}
\end{equation} 

Introducing (\ref{bound2}) into (\ref{bound1}) then yields a bound on the
application between successive moments after harvest.
\begin{equation}\label{bound3}
G(\tilde{x}((n+1)T_h^+))
\leq G(\tilde{x}(nT_h^+))+\int_{nT_h}^{(n+1)T_h}\left[
  \frac{f(\tilde{x}(s))}{g(\tilde{x}(s))}- \min(0,\tilde y_0) e^{-ds} -y_{ph}(s) \right]ds+\frac{\ln(1-\alpha_x)}{r_g}
\end{equation}

We can now evaluate an upper-bound for $G$ at any time $t\geq 0$. Defining $l$
as the integer part of $\frac{t}{T_h}$, we have:
\[
\begin{array}{lll}
G(\tilde{x}(t))-G(x_0)&\leq& \int_{0}^{t}\left[
  \frac{f(\tilde{x}(s))}{g(\tilde{x}(s))}- \min(0,\tilde y_0) e^{-ds} -y_{ph}(s) \right]ds+l
\frac{\ln(1-\alpha_x)}{r_g}\\[.2cm]
&\leq&\int_{0}^{t}\left[
  S_g- \min(0,\tilde y_0) e^{-ds} -y_{ph}(s) \right]ds+l
\frac{\ln(1-\alpha_x)}{r_g}\\[.2cm]
&=&-\int_{0}^{t}\min(0,\tilde y_0) e^{-ds}ds+\int_{lT_h}^{t}\left[
  S_g-y_{ph}(s) \right]ds+l\int_{0}^{T_h}\left[
  S_g-y_{ph}(s) \right]ds\\[.1cm]
&&+l\frac{\ln(1-\alpha_x)}{r_g}\\[.2cm]
&=&\frac{\min(0,\tilde y_0)}{d}(e^{-dt}-1)+\int_{lT_h}^{t}\left[
  S_g-y_{ph}(s) \right]ds+l\int_{0}^{T_h}\left[
  S_g-y_{ph}(s) \right]ds\\[.1cm]
&&+l
\frac{\ln(1-\alpha_x)}{r_g}
\end{array}
\]
The first two terms are bounded (the first one is obvious and the second one
is upper-bounded by $S_gT_h$). We then have to analyze the third one, which has
been obtained through the periodicity of $y_{ph}(t)$ and the fourth in order to
know if $G(\tilde{x}(t))$ goes to $-\infty$ when $t$ goes to infinity. In fact, it
suffices to have
\[
\int_{0}^{T_h}\left[
  S_g-y_{ph}(s) \right]ds+\frac{\ln(1-\alpha_x)}{r_g}<0
\]
to achieve this. It is more cleanly rewritten in the form
\begin{equation}\label{cond_B11_hi} 
\int_{0}^{T_h} y_{ph}(t)dt>S_gT_h+\frac{\ln(1-\alpha_x)}{r_g}
\end{equation}
In order to obtain (\ref{condTh}), we are now left with the computation of
$\int_{0}^{T_h} y_{ph}(t)dt$, which is detailed in Proposition \ref{int_h} of the
Appendix:
\begin{equation}\label{yp_int_hi}
\int_{0}^{T_h} y_{ph}(t)dt=\frac{\mu T_h}{d}\left(1- \left(\frac{\alpha_y\left(1 -
              e^{-dT_h}\right)}{1 - \left(1 -
              \alpha_y\right)e^{-dT_h}  }\right)\left(\frac{ e^{-\frac{d T_h}{k}}}{k\left(1 - e^{-\frac{dT_h }{k}}\right)  }   \right)\right)
\end{equation}

Introducing (\ref{yp_int_hi}) into (\ref{cond_B11_hi}) then yields
(\ref{condTh}), which shows that this last condition is sufficient for having
$\tilde x$ going to $0$ as $t$ goes to $\infty$. 

Since $\tilde x$ goes to zero, there exists a finite time $t_f$ after which
$h(\tilde x)\leq \frac{d}{2}$ for all times. Therefore, after this time, we have
\[
\dot{\tilde y}=h(\tilde{x})(y_{ph}(t)+\tilde{y})-d\tilde{y}\leq h(\tilde{x})y_{ph}(t)-\frac{d}{2}\tilde{y}
\] 
We have seen that $h(\tilde{x})y_{ph}(t)$ goes to zero as $t$ goes to
infinity; so does also $\tilde y$.

In order to have the global asymptotic stability, we are only left with the
local asymptotic stability to prove. In order to do that, we only have to
consider the discrete system that maps the state at time $nT_h^+$ onto
the state at time $(n+1)T_h^+$ with respect to the linear equation
(\ref{linear}) and the discrete part. After some computations, we obtain:

\begin{equation}
\label{til_ref_hi}
 \left(
	\begin{array}{c}
		\tilde{x}\\
		\tilde{y}
	\end{array}
  \right)
  \left(
	(n + 1)T_h^+
  \right)
  =
  \mathbf{B}
  \left(
	\begin{array}{c}
		\tilde{x}\\
		\tilde{y}
	\end{array}
  \right)
  \left( 
	nT_h^+ 
  \right)
\end{equation}
where 
\begin{equation*}
\mathbf{B} 
= 
\left(
	\begin{array}{cc}
		(1 - \alpha_x )e^{\int_{nT_h}^{(n + 1)T_h}{f'(0) - g'(0)y_{ph}}{d\tau}} & 0\\
		\ddagger & (1 - \alpha_y)e^{-d\int_{nT_h}^{(n +1)T_h} {}{d\tau}}  
	\end{array}
\right)
\end{equation*}

Note that $\ddagger$ is a term that we do not use in our analysis, therefore is not expressed. Indeed, since the matrix is triangular, it is stable if $|B_{11}| < 1$, i.e.

\begin{equation}
  	\label{cond_B11_hi_lin}
	\displaystyle\int_{nT_h}^{(n + 1)T_h}{y_{ph}}{d\tau} > \displaystyle\frac{f'(0)T_h + \ln(1 - \alpha_x)}{g'(0)}
\end{equation}

Similarly to what was done earlier, it can be shown that
(\ref{cond_B11_hi_lin}) is equivalent to (\ref{condThL}), so that the
necessary and sufficient condition for local stability is proven. 

It is directly seen that (\ref{cond_B11_hi_lin}) is satisfied when
(\ref{cond_B11_hi}) is because $\underline{\mu}_h(S,r)$ is increasing in $S$
and $r$ and we have

\begin{equation}\label{prime_sup}
\frac{f'(0)}{g'(0)}=\lim_{x\stackrel{\geq}{\rightarrow} 0}\frac{f(x)}{g(x)} \leq\sup_{x\geq
  0}\frac{f(x)}{g(x)} \mbox{ and } g'(0)=\lim_{x\stackrel{\geq}{\rightarrow} 0}\frac{g(x)}{x}\leq\sup_{x\geq
  0}\frac{g(x)}{x}
\end{equation}

This completes the proof of global stability, since we have shown global
convergence and local stability when (\ref{condTh}) is satisfied.

\end{proof}

\subsubsection*{Releases less frequent than harvests}

If we now consider the case where predators releases take place less often
than harvests, we also obtain global and local stability results based on the
following function
\[
	\underline{\mu}_r(S,r)=
			d\left(
				S + \frac{ln(1 - \alpha_x)}{rT_h}  
			\right)
			\frac{
				1 - (1 - \alpha_y)e^{-dT_h}
			}
			{
				1 - e^{-dT_h}
			}
\]
which is increasing in $S$ and $r$ since the last fraction is positive and
$\alpha_x\leq 1$.
 \begin{theorem}\label{theoremTr}
When $T_r=kT_h$ with $k\,\in\mathbb{N^*}$, the solution $(x(t),y(t))=(0,y_{pr}(t))$
of (\ref{model}) is LAS iff
\begin{equation}\label{condTrL}
  \mu > \underline{\mu}_r\left(\frac{f'(0)}{g'(0)},g'(0)\right)  
\end{equation}
and is GAS if
\begin{equation}
\label{condTr}
  \mu > \underline{\mu}_r\left(\sup_{x\geq
  0}\frac{f(x)}{g(x)},\sup_{x\geq
  0}\frac{g(x)}{x}   \right)  
\end{equation}
\end{theorem}

\begin{proof}
This proof does not depart very much from the one of Theorem
\ref{theoremTh}. The only difference is that the reference period is now
$T_r$. We use the same function $G(\tilde{x})$ as in (\ref{Gdef}) and an analysis
identical to the one of the previous theorem leads to
\[\begin{array}{ll}
G(\tilde{x}(mT_r+(l+1)T_h^+))\leq& G(\tilde{x}(mT_r+lT_h^+))\\[.1cm]&+\int_{mT_r+lT_h}^{mT_r+(l+1)T_h }\left[
  \frac{f(\tilde{x}(s))}{g(\tilde{x}(s))}- \min(0,\tilde y_0) e^{-ds} -y_{pr}(s) \right]ds+
\frac{\ln(1-\alpha_x)}{r_g}\end{array}
\]
which is exactly (\ref{bound3}) since it is depicting the
behaviour of the model between two harvesting instants.

Extending this to the whole $T_r$ interval, we obtain
\[\begin{array}{lll}
G(\tilde{x}((m+1)T_r^+))&\leq& G(\tilde{x}(mT_r^+))+\int_{mT_r}^{(m+1)T_r }\left[
  \frac{f(\tilde{x}(s))}{g(\tilde{x}(s))}- \min(0,\tilde y_0) e^{-ds} -y_{pr}(s) \right]ds\\[.1cm]&&+
k\frac{\ln(1-\alpha_x)}{r_g}\end{array}
\]
We now see that this expression is identical to (\ref{bound3})
with the exception of the presence of a $k$ factor and the expression of
$y_{pr}(t)$, which comes from (\ref{theo_periodicsol_lo}) instead of (\ref{theo_periodicsol_hi}). 

Condition (\ref{cond_B11_hi}) then becomes
\begin{equation}\label{cond_B11_ri} 
\int_{0}^{T_r} y_{pr}(t)dt>S_gT_r+k\frac{\ln(1-\alpha_x)}{r_g}
\end{equation}
and the computation of 
$\int_{0}^{T_r} y_{pr}(t)dt$ with $y_{pr}(t)$ as in (\ref{theo_periodicsol_lo}) yields:
\begin{equation}\label{yp_int_ri}
\int_{0}^{T_r} y_{pr}(t)dt= \frac{\mu T_r}{d}\frac{1-e^{-dT_h}}{1-(1-\alpha_y)e^{-dT_h}}
\end{equation}
This leads to conditon (\ref{condTr}) (the computation of (\ref{yp_int_ri}) is
detailed in Proposition \ref{int_r} in the appendix). Global convergence of $(\tilde x,\tilde y)$
to $(0,0)$ is then concluded by
using the same argument as in the proof of Theorem \ref{theoremTh} to show the
convergence of $\tilde y$ to $0$. 

The local stability condition (\ref{condTrL}) then directly arises from the
analysis of the stability of the discrete linearized system that maps $\tilde y(mT_r^+)$
onto $\tilde y((m+1)T_r^+)$. 

\begin{equation}
\left(
	\begin{array}{c}
		\tilde{x}\\
		\tilde{y}
	\end{array}
\right)
\left(
	(m + 1)T_r^+
\right)
=
\mathbf{B}
\left(
	\begin{array}{c}
		\tilde{x}\\
		\tilde{y}
	\end{array}
\right)
\left(
	mT_r^+
\right)
\end{equation}
where
\begin{equation}
\label{B_lo} 
\mathbf{B} = 
\left(
	\begin{array}{cc}
		(1 - \alpha_x)^ke^{\int_{mT_r}^{(m + 1)T_r}{f'(0) - g'(0)y_{pr}}{d\tau}} & 0\\
		e^{\int_{mT_r}^{(m + 1)T_r}{h'(0)y_{pr}}{d\tau}} & (1 - \alpha_y)^ke^{-d\int_{mT_r}^{(m + 1)T_r}{}{d\tau}}  
	\end{array}
\right)
\end{equation}

Again, since the system matrix is lower triangular, for stability we simply
require that $|B_{11}| < 1$ and $|B_{22}| < 1$. The latter yields a trivial
condition, so we calculate the former:
\begin{equation}
  	\label{cond_B11_ri_lin}
	\displaystyle\int_{mT_r}^{(m + 1)T_r}{y_{pr}}{d\tau} > \displaystyle\frac{f'(0)T_r + k\ln(1 - \alpha_x)}{g'(0)}
\end{equation}
which directly leads to condition (\ref{condTrL}) for local stability and the
proof of global stability is also complete because $\underline{\mu}_r$ is
increasing in $S$ and $r$ and (\ref{prime_sup}) is still satisfied.
\end{proof}
\subsubsection*{Comment}
  As we have seen, when the condition (\ref{condTh}) or (\ref{condTr}) is satisfied, the extinction of the pests is
GAS. When the local condition (\ref{condThL}) or (\ref{condTrL}) is not verified, the extinction of the pests is not stable and a
bifurcation analysis similar to what is done in \cite{LAKARI00,LIUZHACHE05a} would
show the presence of a limit cycle when $\mu$ is close to the limit. When
$\mu$ satisfies condition (\ref{condThL}) or (\ref{condTrL}) only, the pests 
extinction is locally stable and we
cannot rule out that it is globally stable (since our global condition is only
sufficient). Such a budget has the advantage of being  smaller than the one
that guarantees global stability. It allows for good control of limited pest
invasions; however the culture is at risk of being destroyed by a large pest
outbreak.

Since, in both cases, the conditions for local and global stability are identical up to two
different parameters, any analysis of the consequences of one of those
conditions will immediately translate to the other. The interpretation of conditions (\ref{condThL})-(\ref{condTh}) and
(\ref{condTrL})-(\ref{condTr}) will be given in the next section.

\section{Interpretation of results}

It is easy to see that $\underline{\mu}_r$ is independent of $T_r$. The influence of $T_r$ on $\underline{\mu}_h$ is trickier to identify so we shall analyse it mathematically first. We then present graphically the variation of both $\underline{\mu}_r$ and $\underline{\mu}_h$ with respect to $T_r$ for a typical set of parameter values, and attenpt to give a practical interpretation of these results.

\subsection{Mathematical analysis}

We first need to note that when $S + \frac{\ln(1 - \alpha_x)}{rT_h} < 0$, for any of the local or global condition, the condition is trivially verified. Indeed, it implies simply that no biological control is needed for exterminating the pests; in fact, the partial harvesting is effective enough for this purpose (as $\alpha_x$ is large enough). We now evaluate how the release frequency influences the minimal budget when this condition is not trivial.

We have already seen that $\underline{\mu}_r$ is independent of $T_r$. We will now study the latter's influence on $\underline{\mu}_h $.


\begin{theorem}
Let $T_h = kT_r$ where $k \in \mathbb{N^*}$. 

The minimal budget is monotonically decreasing with respect to the release period $T_r$ for non-negative values of $\underline{\mu_h}$, i.e.

\begin{equation}
\displaystyle\frac{\partial{\underline{\mu}_h}}{\partial{T_r}} <  0
\end{equation}. 



\end{theorem}

\begin{proof}
Knowing that $T_r$ is equal to $\frac{T_h}{k} $, it is possible to identify the sign of $\frac{\partial{\underline{\mu}_r}}{\partial{T_r}}$ noting that

\begin{equation*}
\label{dmu_r_dTr}
\begin{array}{ccl}
\displaystyle\frac{\partial{\underline{\mu}_h}}{\partial{T_r}}  =  \displaystyle\frac{\partial{\underline{\mu}_h}}{\partial{k}}\displaystyle\frac{\partial{k}}{\partial{T_r}}
  =  \displaystyle\frac{\partial{\underline{\mu}_h}}{\partial{k}}\left(\displaystyle\frac{-k^2}{T_h}\right)
\end{array}
\end{equation*}
So
\begin{equation}
\label{sgn_dmu_r_dTr_1}
\mbox{sgn}{
	\left(
		\displaystyle\frac{\partial{\underline{\mu}_h}}{\partial{T_r}} 
	\right)
	} 
=
-\mbox{sgn}{
	\left(
		\displaystyle\frac{\partial{\underline{\mu}_h}}{\partial{k}} 
	\right)
	}   
\end{equation}
$\underline{\mu}_h$ is expressed as the product of two distinctive parts, one of which is independent of $k$ and which, for the non-trivial stability condition, is positive, 
$$
 S + \displaystyle\frac{\ln(1 - \alpha_x)d}{rT_h} > 0
$$
where $S$ and $r$ are the parameters required for the local and global conditions, as defined previously.

The second part is viewed as a composite function of $k$ so that (\ref{sgn_dmu_r_dTr_1}) can be evaluated as

\begin{equation}
\label{sgn_dmu_r_dTr_2}
 \mbox{sgn}{
  \left(
   \displaystyle\frac{\partial{\underline{\mu}_h}}{\partial{T_r}} 
  \right)}
=
-\mbox{sgn}{
  \left(
    \displaystyle\frac{\partial{}}{\partial{k}}	
    \left(
      \displaystyle\frac{1}{1 - \left(
				  \frac{\alpha_y(1 - e^{-dT_h})}{1 - (1 - \alpha_y)e^{-dT_h}}
				\right)
				\sigma(k)}
    \right)
  \right)}   
\end{equation}
where $\sigma(k) = \left( \frac{e^{-dT_h/k}}{k\left(1 - e^{-dT_h/k}\right)} \right)$. Then, we get

\begin{equation}
\label{sgn_dmu_r_dTr_3}
\begin{array}{ccl}
\mbox{sgn}{
  \left(
   \displaystyle\frac{\partial{\underline{\mu}_h}}{\partial{T_r}} 
  \right)}
& = &
-\mbox{sgn}{
  \left(
      \displaystyle\frac{\frac{\alpha_y(1 - e^{-dT_h})}{1 - (1 - \alpha_y)e^{-dT_h}}}{\left(
				   1 - \left(
				   \frac{\alpha_y(1 - e^{-dT_h})}{1 - (1 - \alpha_y)e^{-dT_h}}
					\right)
			    \sigma(k)\right)^2
				}
     \displaystyle\frac{\partial{\sigma}}{\partial{k}}	
  \right)} \\
& = & 
-\mbox{sgn}{
  \left(
    \displaystyle\frac{\partial{\sigma}}{\partial{k}}	
  \right)} \\
& = &
-\mbox{sgn}\left({
  \displaystyle\frac{e^{-dT_h/k}}{k^2(1-e^{-dT_h/k})^2}\left(\displaystyle\frac{dT_h}{k} - 1 + e^{-dT_h/k} \right)
}\right)\\\\
& = & 
-\mbox{sgn}{
  \left(
    ke^{-dT_h/k} + dT_h - k
  \right)
}
\end{array}
\end{equation}
Since
\begin{equation*}
\begin{array}{ccl}
\displaystyle\frac{\partial{}}{\partial{k}}(ke^{-dT_h/k} + dT_h - k) 
& = & 	\left(
	  1 + \displaystyle\frac{dT_h}{k}
	\right)e^{-dT_h/k} - 1\\
& \leq & 0 
\end{array}
\end{equation*}
and using l'Hospital's Rule
\begin{equation*}
\begin{array}{ccl}
\displaystyle\lim_{k\to\infty}{\left(ke^{-dT_h/k} + dT_h - k\right)} 
& = &
dT_h + \displaystyle\lim_{k\to\infty}{\left( \displaystyle\frac{e^{-dT_h/k} - 1}{\frac{1}{k}} \right)}\\\\
& = & 
 dT_h + \displaystyle\lim_{k\to \infty}{\left( \displaystyle\frac{\frac{dT_h}{k^2}e^{-dT_h/k}}{-\frac{1}{k^2}} \right)}\\
& = &
0
\end{array}
\end{equation*}
we deduce that $\mbox{sgn}{
  \left(
    ke^{-dT_h/k} + dT_h - k
  \right)
} > 0$
Therefore,
\begin{equation}
\label{sgn_dmu_r_dTr_4}
  \begin{array}{ccl}
    \mbox{sgn}{
    \left(
      \displaystyle\frac{\partial{\underline{\mu}_h}}{\partial{T_r}}  
    \right)}
    < 0
  \end{array}
\end{equation}

\end{proof}

We can deduce that we hit the smallest minimal value for the budget for the largest possible $T_r$ in this case that corresponds to when $k = 1$. This happens when the release frequency equals the partial harvest frequency.

\subsection{Discussion}
\begin{figure}[ht!]
\begin{center}
\includegraphics[width=0.75\textwidth]{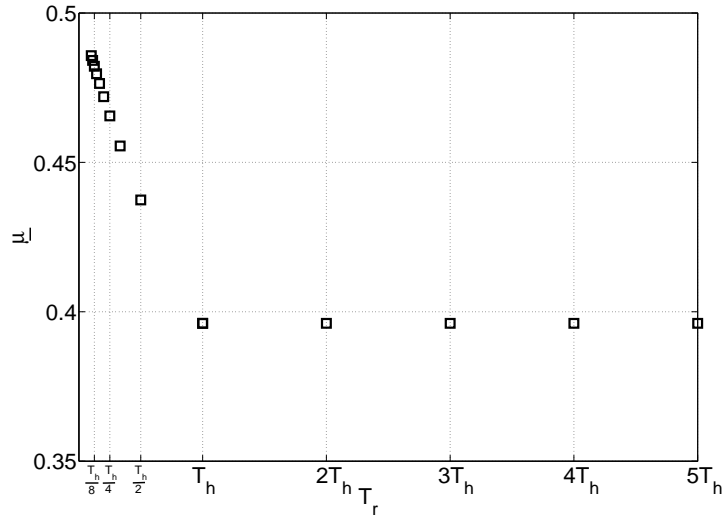}
\caption{Variation of the minimal number of predators required per budget year $\underline{\mu}$ as a function of release to harvest period ratio. Parameters are given the values (in arbitrary units): $\alpha_x = \alpha_y =  0.5$, $d = 1$, and the rate of growth $f'(0)$, functional response $g'(0)$ and numerical response $h'(0)$ with respect to time when the ecosystem is pest-free, i.e. $x_p(t) = 0$, are all equal to 1. }\label{fig_muinf_Tr}
\end{center}
\end{figure}

 Figure \ref{fig_muinf_Tr} represents the analytical results obtained in the previous sections for a chosen set of parameters. The plot includes the two studied cases: either one of the partial harvest and the release period is an integer multiple of the other. 

Under this set of possible scenarios, increasing the frequency of release beyond the frequency of harvest requires that the total number of predators to invest in be higher than that when releases take place less or as often as partial harvests. In the latter case, to ensure pest eradication, the total budget of predators to invest in is fixed, independently of the release period. 

These results imply that it is clearly less costly to protect a greenhouse culture for lower frequencies of release. Of additional economic interest, in this case, the biological treatment is always combined with partial harvesting, so that there is little or no extra cost linked to the presence of workers on-site.  However, we recall that \cite{MAIGRO06} previously demonstrated that the higher the release frequency, the smaller the worst-case damages. Combining the results from both studies seems to indicate that the most profitable release strategy among the possibilities that have been considered is the one where releases are synchronized with the partial harvests.

\section{Conclusion}

The results obtained in this paper for the stability of the system are yet another confirmation that inundative control can be an effective means of suppressing low pest invasions in a greenhouse. This requires that a sufficient number of predators are introduced in the system as in, for instance, \cite{FENeA01,JACCROFEN01}. 

Our study aimed to provide a control strategy in the protection of continuously grown crops that are partially harvested on a regular basis. We demonstrated that partial harvesting had a non-negligible effect on biological control and needed to be taken into account when devising a control strategy in the case of such crops.

We thus investigated the combined effects of releases and partial harvests in terms of the relative frequencies of their implementation. We considered the case where these two events occurred at periods such that one was the integer multiple of the other, and with the two events coinciding over the longer period. In particular, we found when releases were as frequent as or less than the partial harvests, the minimal budget did not depend on the period of release but instead on the harvest parameters, the growth function of the pest population, the mortality of the predators and the functional response. When releases were more frequent than the partial harvests however, the minimal budget value increased with the increasing frequency of the releases, exceeding the constant value obtained for the less frequent case. Combined with the findings of \cite{MAIGRO06} which pointed out that higher release frequencies led to the optimal control policy, we concluded that for the set of possibilities that was studied, the current best strategy is when release and harvest frequencies are equal.

This approach has, however, its shortcomings. Since the integer multiple factor is key to calculating the minimal budget which would satisfy the stability conditions, it is not yet generalised to other scenarios where neither period is the integer multiple of the other. This would happen for instance at other rational non-integer ratios as irrational ones. It is highly likely that these intermediate ratios might induce other dynamics in the system. Whether they might stabilise it given even lower minimal budget values or favour chaos remains to be seen. Moreover it would be interesting to extend the results to the case where the two controls never coincide in spite of following a periodic pattern. This would be in the line of the work, for instance, of \cite{LIUCHEZHA05}, where pesticide spraying - which is analogous to harvests - and releases are not synchronised. 

Nevertheless, we consider that our simplification already has its practical economical advantage. Indeed, coinciding periods imply little or no additional costs incurred in terms of labour: the task of predator release can be assigned to workers in charge of partial harvesting. Field-testing is now the next step required to validate the results of this paper.

\newpage
\section*{Appendix}
\begin{proposition}\label{int_h} Let Hypothesis \ref{hyp1} hold, then
$$
\displaystyle\int_{0}^{T_h} y_{ph}(t)dt=\frac{\mu T_h}{d}\left(1- \left(\frac{\alpha_y\left(1 -
              e^{-dT_h}\right)}{1 - \left(1 -
              \alpha_y\right)e^{-dT_h}  }\right)\left(\frac{ e^{-\frac{d T_h}{k}}}{k\left(1 - e^{-\frac{dT_h }{k}}\right)  }   \right)\right)
$$
\end{proposition}
\begin{proof}
In order to compute the integral, we describe $y_{ph}(t)$ as
$y_{ph}(iT_r^+)e^{-(t-iT_r)}$ in each time interval $[iT_r,(i+1)T_r]$, with
$y_{ph}(iT_r^+)$ given by (\ref{y_Tr}) when $y(nT_h^+)=y^*$. This yields:
\[
\begin{array}{ll}
\int_{0}^{T_h} y_{ph}(t)dt&= \displaystyle\sum_{i=0}^{k-1} y_{ph}(iT_r^+) \int_{iT_r}^{(i + 1)T_r}e^{-d(t-iT_r)}dt\\
&=\displaystyle\sum_{i=0}^{k-1} \left(y^*e^{-idT_r}+\mu
T_r\displaystyle\sum_{j=0}^{i-1}e^{-jdT_r}\right)\int_{0}^{T_r}e^{-dt}dt\\
&=\left(y^*\frac{1-e^{-kdT_r}}{1-e^{-dT_r}}+\mu
  T_r\sum_{i=0}^{k-1}\frac{1-e^{-idT_r}}{1-e^{-dT_r}}\right)\frac{1-e^{-dT_r}}{d}\\
&=\frac{y^*}{d}(1-e^{-kdT_r})+\frac{\mu T_r}{d}\left(k-\frac{1-e^{-kdT_r}}{1-e^{-dT_r}} \right)\\
&=\frac{	
			\left(
				\frac{1 - e^{-dT_h}}{1 - e^{-dT_r}} 
			\right) (1 - \alpha_y) + \alpha_y }{1 -(1 - \alpha_y)e^{-dT_h} }
		\frac{\mu T_r}{d}(1-e^{-kdT_r})+\frac{\mu
                  T_r}{d}\left(k-\frac{1-e^{-kdT_r}}{1-e^{-dT_r}} \right)\\
&=\mu\frac{T_h}{dk}\frac{\left(\left( 1 - e^{-dT_h}\right)(1 -
    \alpha_y)+\alpha_y\left(1 - e^{-\frac{dT_h}{k}} \right)\right)  \left( 1 -
    e^{-dT_h}\right)}{\left(1 - e^{-\frac{dT_h}{k}} \right)\left(1 -(1 - \alpha_y)e^{-dT_h}
  \right) } \\
& \hspace{4cm} + \mu\frac{T_h}{dk}\frac{\left(k \left(1 - e^{-\frac{dT_h}{k}} \right)-  \left( 1
      - e^{-dT_h}\right)\right) \left(1 -(1 - \alpha_y)e^{-dT_h} \right)}
{\left(1 - e^{-\frac{dT_h}{k}} \right)\left(1 -(1 - \alpha_y)e^{-dT_h}
  \right)            }\\
&=\mu\frac{T_h}{dk}\frac{\left(\left( 1 - (1-\alpha_y)e^{-dT_h}\right)-\alpha_y e^{-\frac{dT_h}{k}}\right)  \left( 1 -
    e^{-dT_h}\right)}{\left(1 - e^{-\frac{dT_h}{k}} \right)\left(1 -(1 - \alpha_y)e^{-dT_h}
  \right) }\\
&\hspace{4cm}+ \mu\frac{T_h}{dk}\frac{\left(k \left(1 - e^{-\frac{dT_h}{k}} \right)-  \left( 1
      - e^{-dT_h}\right)\right) \left(1 -(1 - \alpha_y)e^{-dT_h} \right)
}{\left(1 - e^{-\frac{dT_h}{k}} \right)\left(1 -(1 - \alpha_y)e^{-dT_h}
  \right) }\\
&=\mu\frac{T_h}{dk}\frac{\left(-\alpha_y e^{-\frac{dT_h}{k}}\right)  \left( 1 -
    e^{-dT_h}\right)+\left(k \left(1 - e^{-\frac{dT_h}{k}} \right)\right) \left(1 -(1 - \alpha_y)e^{-dT_h} \right)
}{\left(1 - e^{-\frac{dT_h}{k}} \right)\left(1 -(1 - \alpha_y)e^{-dT_h}
  \right) }\\
&=\frac{\mu T_h}{d}\left(1- \left(\frac{\alpha_y\left(1 -
              e^{-dT_h}\right)}{1 - \left(1 -
              \alpha_y\right)e^{-dT_h}  }\right)\left(\frac{ e^{-\frac{d T_h}{k}}}{k\left(1 - e^{-\frac{dT_h }{k}}\right)  }   \right)\right)
\end{array}
\]
\end{proof}
\begin{proposition}\label{int_r} 
Let Hypothesis \ref{hyp1} hold, then 
\[
\int_{0}^{T_r} y_{pr}(t)dt= \frac{\mu
  T_r}{d}\frac{1-e^{-dT_h}}{1-(1-\alpha_y)e^{-dT_h}}
\]
\end{proposition}
\begin{proof}
In order to compute the integral, we describe $y_{pr}(t)$ as
$y_{pr}(iT_h^+)e^{-(t-iT_h)}$ in each time interval $[iT_h,(i+1)T_h]$, with
$y_{pr}(iT_h^+)$ given by (\ref{y_Th}) when $y(mT_r^+)=y^*$. This yields:
\[
\begin{array}{lll}
\int_{0}^{T_r} y_{pr}(t)dt&=& \displaystyle\sum_{i=0}^{k-1} y_{pr}(iT_h^+)
\int_{iT_h}^{(i + 1)T_h}e^{-d(t-iT_h)}dt\\
&=&\displaystyle\sum_{i=0}^{k-1} y^* e^{idT_h}(1-\alpha_y)^i
\int_{0}^{T_h}e^{-dt}dt\\
&=&y^*\frac{(1-e^{-dT_h})}{d}\displaystyle\sum_{i=0}^{k-1} y^* e^{idT_h}(1-\alpha_y)^i\\
&=&\frac{\mu T_r}{1 - (1 -
  \alpha_y)^ke^{-dT_r}}\frac{(1-e^{-dT_h})}{d}\frac{1-(1-\alpha_y)^ke^{-kdT_h}}{1-(1-\alpha_y)e^{-dTh}}\\
&=&\frac{\mu T_r(1-e^{-dT_h})}{d(1-(1-\alpha_y)e^{-dTh})}
\end{array}
\]
\end{proof}

\newpage
\tableofcontents

\bibliographystyle{plain}
\bibliography{PhD_bib}

\end{document}